\documentclass{amsart}
\usepackage[all]{xy}
\usepackage{amsmath}
\usepackage{amsfonts,amssymb}
\usepackage{latexsym}
\RequirePackage[T1]{fontenc}

\newtheorem{dfn}{Definition}[section]
\newtheorem{prop}[dfn]{Proposition}

\newtheorem{cjt}[dfn]{Conjecture}
\newtheorem{thm}[dfn]{Theorem}
\newtheorem{lem}[dfn]{Lemma}

\newtheorem{exs}[dfn]{Examples}

\newcommand{\End}[1]{\textrm{End}(#1)}

\newcommand{\R}{\mathbb{R}}

\newcommand{\N}{\mathbb{N}}

\newcommand{\al}{\mathfrak g}
\newcommand{\gun}{\mathfrak g_1}
\newcommand{\gde}{\mathfrak g_2}
\newcommand{\gl}{\mathfrak{gl}}

\author{Damien Calaque}

\title{Formality for Lie algebroids}

\begin{document}

\maketitle

\noindent{\small IRMA, 7 rue René Descartes, F-67084 Strasbourg, France \\
\emph{E-mail address}: {\bf calaque@math.u-strasbg.fr}}

\begin{abstract}
Using Dolgushev's generalization of Fedosov's method for deformation quantization, we give a positive answer to a question 
of P.~Xu: can one prove a formality theorem for Lie algebroids ? 
As a direct application of this result, we obtain that any triangular Lie bialgebroid is quantizable. 
\end{abstract}

\section*{Introduction}

The main goal of this paper is to formulate precisely a `Kontsevich-like' (see \cite{K}) formality theorem for Lie 
algebroids and then prove it. This problem has been proposed by Ping Xu at the end of \cite{X2} (question 2). To do it, 
we use a slightly modified version of Dolgushev's equivariant globalisation of Kontsevich's formality \cite{D} (in his 
paper Dolgushev generalises Fedosov's geometric construction of star-products \cite{F} to the case of a general manifold). 

We would like to mention that most of the proofs given in this paper are similar to those of \cite{D} and apologize for this 
repetition. 

The paper is organized as follows. 

Section 1 is devoted to the presentation of our main results. We first recall some basic facts about differential geometry 
for Lie algebroids (see \cite{CW,M} for details). Then we prove a Hochschild-Kostant-Rosenberg theorem for Lie algebroids 
and state our formality theorem. We finally explain why this result implies that any triangular Lie bialgebroid is 
quantizable (which has been proved in the case of regular ones by Ping Xu \cite{X2}, see also \cite{NT}). 

In section 2, we construct resolutions of the desired DGLA using a torsion free Lie algebroid connexion. It is the more 
technical part of the paper. 

We end the proof of our main theorem in section 3: after twisting a fiberwize quasi-isomorphism (given by \cite{K}) we use 
the resolutions of the previous section to contract it to the desired one. We also prove an equivariant version of our 
result. 

We recall in an appendix some facts about $L_\infty$-algebras, Hopf algebroids, and Lie algebroid connections. 

Throughout the paper the Einstein convention for summation over repeated indices is assumed. \\

\emph{Acknowledgements. }
I am grateful to my advisor, B.~Enriquez, who has accepted to lead my research and 
read carefully this paper. I am also greatly indebted to G.~Halbout for teaching me the ideas of \cite{D,F,K}. Discussions 
on `thing-oids' with P.~Xu in Normandie were very enlightening, I express to him my sincere thanks. I also thank 
V.~Dolgushev for his warned comments. 

\section{Main results}

\subsection{Preliminaries: differential geometry for Lie algebroids}\label{sec:pre}

\begin{dfn}[\cite{M}]\emph{
A \emph{Lie algebroid} is a vector bundle $E$ over a manifold $X$ equipped with a Lie bracket $[,]_E$ on sections 
$\Gamma(X,E)$ and a bundle map $\rho:E\to TX$ called the \emph{anchor} such that: 
\begin{enumerate}
\item The induced map $\rho:\Gamma(X,E)\to\Gamma(X,TX)$ is a Lie algebra morphism. 
\item For any $f\in C^\infty(X)$, $v,w\in\Gamma(X,E)$, 
$$[v,fw]_E=f[v,w]_E+(\rho(v)\cdot f)w \qquad (\textrm{Leibniz identity}) $$
\end{enumerate}
}\end{dfn}

Basic objects in differential geometry are tensors. So it is natural to consider their algebroids analogues which we call 
\emph{$E$-tensors}: for $k,l\geq0$, an $E$-$(k,l)$-tensor is a section of the bundle $(\otimes^kE)\otimes(\otimes^lE^*)$. 
In a local base $(e_1,\dots,e_r)$ of $E$ with duale base $(\xi^1,\dots,\xi^r)$ of $E^*$, such an $E$-tensor $T$ can be 
written 
$$T(x)=T_{j_1\dots j_l}^{i_1\dots i_k}(x)e_{i_1}\otimes\cdots\otimes e_{i_k}\otimes\xi^{j_1}\otimes\cdots\otimes\xi^{j_l}$$
Indices $i_1,\dots,i_k$ and $j_1,\dots,j_l$ are respectively called contravariant and covariant. 

As in usual differential geometry, one can consider the graded commutative algebra of \emph{$E$-differential forms} 
$^E\Omega(X):=\Gamma(X,\wedge E^*)$, which is endowed with a square zero super-derivation 
$d_E:~\!\!\!^E\Omega^*(M)\to~\!\!\!^E\Omega^{*+1}(M)$. In local $E$-coordinates, any $E$-$k$-form $\omega$ can be 
written 
$$\omega(x)=\omega_{i_1\dots i_k}(x)\xi^{i_1}\wedge\cdots\wedge\xi^{i_k}$$
where $\omega_{i_1\dots i_k}$ are coefficients of a covariant $E$-tensor antisymmetric in indices $i_1,\dots,i_k$, and 
$$d_E=\xi^i\rho(e_i)-\frac12\xi^i\wedge\xi^jc_{ij}^k(x)\frac\partial{\partial\xi^k}$$
where $[e_i,e_j]_E=c_{ij}^k(x)e_k$. 

In the same way, one can define the differential graded Lie algebra (DGLA for short) of \emph{$E$-polyvector fields} 
$$T_{poly}E=\bigoplus_{k\geq-1}T^k_{poly}E=\bigoplus_{k\geq-1}\Gamma(X,\wedge^{k+1}E)$$
endowed with the zero differential and the Lie super-bracket of degree zero which extend uniquely $[,]_E$ as follow: for 
$u\in T_{poly}^kE$, $v\in T_{poly}^lE$ and $w\in T_{poly}E$, 
$$[u,v\wedge w]_E=[u,v]_E\wedge w+(-1)^{k(l+1)}v\wedge[u,w]_E$$
As in the case of $E$-forms, any $E$-$k$-vector field $v$ can be written locally
$$v(x)=v^{i_1\dots i_k}(x)e_{i_1}\wedge\cdots\wedge e_{i_k}$$
with $v^{i_1\dots i_k}$ coefficients of a contravariant $E$-tensor antisymmetric in indices $i_1,\dots,i_k$. \\

The usual algebra of differential operators can be viewed as a kind of universal enveloping algebra of the Lie algebra of 
vector fields. Let us define in a similar way the algebra of \emph{$E$-differential operators} as the quotient of the 
graded algebra freely generated\footnote{We consider here a completed tensor product: infinite sums which are finite on any 
compact are allowed. } by $C^\infty(X)$ (of degree $0$) and $\Gamma(X,E)$ (of degree $1$) by relations 
{\footnotesize
\begin{eqnarray*}
& f\otimes g-fg & f,g\in C^\infty(X) \\
& f\otimes v-fv & f\in C^\infty(X),v\in\Gamma(X,E) \\
& v\otimes f-f\otimes v-\rho(v)\cdot f & f\in C^\infty(X),v\in\Gamma(X,E) \\
& v_1\otimes v_2-v_2\otimes v_1-[v_1,v_2]_E & v_i\in\Gamma(X,E)
\end{eqnarray*}}
$UE$ carries a natural Hopf algebroid structure (see appendix B) with base algebra $C^\infty(X)$, source and target maps 
$s=t:C^\infty(X)\to UE$ the natural inclusion, coproduct $\Delta:UE\to UE\widetilde\otimes UE$ 
(where $\widetilde\otimes$ denotes $\otimes_{C^\infty(X)}$) which extends 
\begin{equation}\label{eq:coprod}
\begin{array}{l}
\Delta(f)=f\widetilde\otimes1=1\widetilde\otimes f, \forall f\in C^\infty(X) \\
\Delta(v)=v\widetilde\otimes1+1\widetilde\otimes v, \forall v\in\Gamma(X,E)
\end{array}
\end{equation}
and counit $\epsilon:UE\to C^\infty(X)$ which extends 
$$\epsilon(f)=f, \forall f\in C^\infty(X)~~\textrm{and}~~\epsilon(v)=0, \forall v\in\Gamma(X,E)$$
This allows us to define a Lie super-bracket on the graded vector space 
$$D_{poly}E=\bigoplus_{k\geq-1}D_{poly}^kE=\bigoplus_{k\geq-1}UE^{\widetilde\otimes k+1}$$
of $E$-polydifferential operators in a way similar to appendix \ref{sec:A2}: for homogeneous elements 
$P_i\in D^{k_i}_{poly}(X)~(i=1,2)$, $[P_1,P_2]=P_1\bullet P_2-(-1)^{k_1k_2}P_2\bullet P_1$, where
$$P_1\bullet P_2=\sum_{i=0}^{k_1}(-1)^{ik_2}
\textrm{id}^{\widetilde\otimes i}\widetilde\otimes\Delta^{k_2+1}\widetilde\otimes\textrm{id}^{\widetilde\otimes k_1-i}(P_1)
\cdot(1^{\widetilde\otimes i}\widetilde\otimes P_2\widetilde\otimes1^{\widetilde\otimes k_1-i})$$
Remark that $m_0=1\widetilde\otimes1\in D_{poly}E$ is such that $[m_0,m_0]=0$, thus $(\partial=[m_0,\cdot],[,])$ defines a 
DGLA structure on $D_{poly}E$. By an easy calculation, one can observe that $\partial$ is simply the Hochschild coboundary 
operator (up to a sign) for the complex $\widetilde\otimes^*UE$. 

\subsection{Formality for Lie algebroids}

Let $E\to X$ be a Lie algebroid. \\

First, in the spirit of the Hochschild-Kostant-Rosenberg theorem we are going to prove that 
$H^*(D_{poly}E,\partial)\cong T_{poly}^*E$ (which was first proved in \cite{V} for $E=TX$). 
\begin{thm}\label{thm:hkr}
Define the map $U_{hkr}:(T_{poly}E,0)\to(D_{poly}E,\partial)$ by 
$$U_{hkr}(v_0\wedge\cdots\wedge v_n)=\frac1{(n+1)!}\sum_{\sigma\in S_{n+1}}\epsilon(\sigma)v_{\sigma_0}\otimes\cdots\otimes v_{\sigma_n}$$
if $n\geq0$ and $v_i\in\Gamma(X,E)$, and $U_{hkr}(f)=f$ if $f\in C^\infty(X)$. \\
It is a quasi-isomorphism of complexes (i.e., it is a morphism of complexes which induces an isomorphism in cohomology). 
\end{thm}
\begin{proof}
First, one can immediately check that the image of $U_{hkr}$ is annihilated by $\partial$, i.e.~that it is a morphism of 
complexes. 

Now remark that the complex $D_{poly}E$ is filtered by the total degree of polydifferential operators. $T_{poly}E$ carries 
also a natural filtration (which is in fact a gradation), namely by degree of polyvector fields. Then $U_{hkr}$ is 
compatible with filtrations. Thus we have to prove that $Gr(U_{hkr}):Gr(T_{poly}E)\to Gr(D_{poly}E)$ is a quasi-isomorphism 
of complexes. In $Gr(D_{poly}E)$ all components are sections of some vector bundle on $X$ and $\partial$ is 
$C^\infty(X)$-linear (the same is obviously true for $T_{poly}E$), therefore we have to show that $Gr(U_{hkr})$ is a 
quasi-isomorphism fiberwise. 

Fix $x\in X$ and consider the vector space $V=E_x$. One has
$$Gr(D_{poly}E)_x=\bigoplus_{n\geq0}S(V)^{\otimes n}$$
but it is better to indentify $S(V)$ with the cofree cocommutative coalgebra with counit $\mathcal C:=C(V)\oplus(\R1)^*$. 
As above the differential can be expressed in terms of the cocommutative coproduct $\Delta$; namely 
$$(-1)^{n-1}\partial=1^*\otimes\textrm{id}^{\otimes n}-\sum_{i=1}^{n-1}(-1)^i\textrm{id}\otimes\cdots\otimes\Delta_i
\otimes\cdots\otimes\textrm{id}+(-1)^{n-1}\textrm{id}^{\otimes n}\otimes1^*$$
Now let us recall a standard result in homological algebra: 
\begin{lem}
Let $\mathcal C$ be the cofree cocommutative coalebra with counit cogenerated by a vector space $V$. Then the natural 
homomorphism of complexes $(\wedge^*V,0)\to(\otimes^*\mathcal C,\partial)$ is a quasi-isomorphism. 
\end{lem}
Apply this lemma in the case when $V=E_x$ and remark that $Gr(T_{poly}E)_x=(T_{poly}E)_x=\wedge^*V$. The theorem is proved. 
\end{proof}
Now we claim that $D_{poly}E$ is formal: 
\begin{thm}[Formality]\label{thm:formal}
There exists a quasi-isomorphism $U_E$ of DGLA from $(T_{poly}E,0,[,]_E)$ to $(D_{poly}E,\partial,[,])$. 
\end{thm}
When $E=TX$ this is the formality theorem for manifolds presented in \cite[section 4.6]{K}. More 
generally, if the anchor of $E$ is injective then $E\subset TX$ is an integrable distribution (i.e., $X$ is foliated) 
and thus one obtains a formality theorem for leafwise polydifferential operators. 

\subsection{Quantization of triangular Lie bialgebroids}

Let $E\to X$ be a Lie algebroid with bracket $[,]_E$ and anchor $\rho$. 

Let $H=UE$, $R=C^\infty(X)$, $\Delta$ defined by (\ref{eq:coprod}), $s=t:R\to UE$ be the natural embedding and 
$\varepsilon:UE\to R$ extending
$$\begin{array}{l}
\varepsilon(f)=f, \forall f\in R=C^\infty(X) \\
\varepsilon(v)=0, \forall v\in H=UE
\end{array}$$
It is a Hopf algebroid (see appendix \ref{sec:B}). Moreover one can obviously extend the anchor to a map 
$\rho:UE\to U(TX)\subset\End R$. It defines an anchor for the Hopf algebroid $H$. 

In \cite{X2} Ping Xu observes that any Hopf algebroid deformation of $UE$ endows $E$ with a Lie bialgebroid structure (by 
taking the semi-classical limit). Recall the 
\begin{dfn}[\cite{MX}]\emph{
A \emph{Lie bialgebroid} is a Lie algebroid $E\to X$ whose dual bundle $E^*\to X$ is also a Lie algebroid and such that 
the differential $d_{E^*}$ on $\Gamma(X,\wedge^*E)$ is a derivation of the super-bracket $[,]_E$; namely 
$$\forall v,w\in\Gamma(X,E), d_{E^*}[v,w]_E=[d_{E^*}v,w]_E+[v,d_{E^*}w]_E$$
A Lie bialgebroid $E$ is called \emph{triangular} if $d_{E^*}=[\Lambda,\cdot]_E$ for a given 
$\Lambda\in\Gamma(X,\wedge^2E)$ satisfying $[\Lambda,\Lambda]_E=0$. 
}\end{dfn}
Reciprocally we say that a Lie bialgebroid is \emph{quantizable} if there exists a deformation of $UE$ whose 
semi-classical limit is precisely the starting bialgebroid structure. 
\begin{cjt}
Any Lie bialgebroid is quantizable. 
\end{cjt}
Following \cite{Dr}, Xu shows in \cite{X2} that to quantize a triangular Lie bialgebroid it is sufficient to find a twistor 
(see appendix \ref{sec:B}) $J\in(UE\otimes_RUE)[[\hbar]]$ such that 
$\frac{J-J^{\textrm{op}}}\hbar=\Lambda\textrm{ mod }\hbar$ and consider 
$(UE[[\hbar]],R_J,\Delta_J,s_J,t_J,\varepsilon)$. We construct such a $J$ with the help of our 
formality theorem \ref{thm:formal}. 
\begin{thm}
Any triangular Lie bialgebroid is quantizable. 
\end{thm}
\begin{proof}
Let us define 
$$J=m_0+\sum_{n\geq1}\frac{\hbar^n}{n!}U_E^{[n]}(\Lambda,\dots,\Lambda)$$
Now since $U$ is a $L_\infty$-morphism $\alpha=(J-m_0)\in\hbar(D^1_{poly}E)[[\hbar]]$ is a Maurer-Cartan element, 
$\partial\alpha+\frac12[\alpha,\alpha]=0$. It means that 
\begin{eqnarray*}
0 & = & [m_0,J]-[m_0,m_0]+\frac12([J,J]-[m_0,J]-[J,m_0]+[m_0,m_0]) \\
& = & [m_0,J]+\frac12([J,J]-2[m_0,J])=\frac12[J,J]
\end{eqnarray*}
Then remark that $J^{12,3}J^{1,2}-J^{1,23}J^{2,3}=\frac12[J,J]$. 

Finally, since $U^{[1]}$ is a quasi-isomorphism of complexes we have 
$\Lambda={\rm Alt}(U^{[1]}(\Lambda))=\frac{J-J^{op}}{\hbar}+O(\hbar)$. 
\end{proof}

\section{Dolgushev-Fedosov resolutions of $T_{poly}E$ and $D_{poly}E$}\label{sec:2}

Let $E\to X$ be a Lie algebroid with bracket $[,]_E$ and anchor $\rho$. 

\subsection{The Weyl bundle and related bundles}

Consider the bundle of algebras $\mathcal W=\hat S(E^*)$, whose sections are functions on $E$ formal in the fibers. 
Any section $s\in\Gamma(X,\mathcal W)$ can be written locally
$$s=s(x,y)=\sum_{l=0}^\infty s_{i_1\dots i_l}(x)y^{i_1}\cdots y^{i_l}$$
where $y^i$ are formal coordinates on the fibers of $E$ and $s_{i_1\dots i_l}$ are coefficients of a covariant $E$-tensor 
symmetric in indices $i_1,\dots,i_l$. 

In the same way, one can define the bundle $\mathcal T=\bigoplus_{k\geq-1}\mathcal T^k$ of formal fiberwise polyvector 
fields on $E$; $\mathcal T^k=\mathcal W\otimes\wedge^{k+1}E$. Any homogeneous section $v\in\Gamma(X,\mathcal T^k)$ is 
locally of the form 
$$v=\sum_{l=0}^\infty v_{i_1\dots i_l}^{j_0\dots j_k}(x)y^{i_1}\cdots y^{i_l}
\frac\partial{\partial y^{j_0}}\wedge\cdots\wedge\frac\partial{\partial y^{j_k}}$$
where $v_{i_1\dots i_l}^{j_0\dots j_k}$ are coefficients of an $E$-tensor symmetric in covariant indices $i_1,\dots,i_l$ 
and antisymmetric in contravariant indices $j_0,\dots,j_k$. 

Finally, we denote by $\mathcal D=\bigoplus_{k\geq-1}\mathcal D^k$ the bundle of formal fiberwise polydifferential 
operators on $E$; $\mathcal D^k=\mathcal W\otimes S(E)^{\otimes k+1}$. Any homogeneous section $P\in\Gamma(X,\mathcal D^k)$ 
is locally of the form 
$$P=\sum_{l=0}^\infty P_{i_1\dots i_l}^{\alpha_0\dots\alpha_k}(x)y^{i_1}\cdots y^{i_l}\frac{\partial^{\vert\alpha_0\vert}}
{\partial y^{\alpha_0}}\otimes\cdots\otimes\frac{\partial^{\vert\alpha_k\vert}}{\partial y^{\alpha_k}}$$
where $\alpha_s$ are multi-indices and $P_{i_1\dots i_l}^{\alpha_0\dots\alpha_k}$ are coefficients of an $E$-tensor 
symmetric in covariant indices $i_1,\dots,i_l$. 

For our purposes, we need to tensor these bundles with the exterior algebra bundle $\wedge E^*$. Namely, we need to 
consider the space $^E\Omega(X,\mathcal B)$ of $E$-differential forms on $X$ with values in $\mathcal B$ (from now, 
$\mathcal B$ will denote either $\mathcal W$, $\mathcal T$ or $\mathcal D$). In this setting, $^E\Omega(X,\mathcal W)$ has 
a natural structure of super-commutative algebra, and $^E\Omega(X,\mathcal T)$ (resp.~$^E\Omega(X,\mathcal D)$) is 
naturally endowed with the DGLA structure induced fiber-by-fiber by the DGLA structure of $T_{poly}(\R^d_{formal})$ 
(resp.~$D_{poly}(\R^d_{formal})$). Let us denote the differential and the Lie super-bracket in $^E\Omega(X,\mathcal D)$ by 
$\partial$ and $[,]_G$ respectively, and the Lie super-bracket in $^E\Omega(X,\mathcal T)$ by $[,]_S$. \\

In what follows we denote the same operations on these three different algebras by the same letters when it does not 
lead to any confusion. \\

The differential $\delta=\xi^i\frac\partial{\partial y^i}:~^E\Omega^*(X,\mathcal W)\to~^E\Omega^{*+1}(X,\mathcal W)$
($\delta^2=0$) can obviously extend to $^E\Omega(X,\mathcal T)$ and $^E\Omega(X,\mathcal D)$. Namely, 
$\delta=[\xi^i\frac\partial{\partial y^i},\cdot]_S$ on $^E\Omega(X,\mathcal T)$ and 
$\delta=[\xi^i\frac\partial{\partial y^i},\cdot]_G$ on $^E\Omega(X,\mathcal D)$. 
By definition, $\delta$ is a derivation of the Lie algebras $^E\Omega(X,\mathcal T)$ and $^E\Omega(X,\mathcal D)$. 
Moreover, $\delta$ and $\partial$ super-commute since the multiplication operator $m$ is $\delta$-closed ($\delta m=0$). 
Consequently, $\delta$ is compatible with DGLA structures on $^E\Omega(X,\mathcal T)$ and $^E\Omega(X,\mathcal D)$. 
\begin{prop}
For all $n>0$, $H^n(^E\Omega(X,\mathcal B),\delta)=0$, and $H^0(^E\Omega(X,\mathcal B),\delta)=F^0\mathcal B$ is the space 
of sections of $\mathcal B$ that are constant in the fibers. 
\end{prop}
\begin{proof}
Let us introduce the operator $\delta^*=y^i\iota(e_i)$ of contraction with the Euler vector field $\Theta_E=y^ie_i$. On 
sections of $^E\Omega^k(X,\mathcal B)$ polynomial of degree $l$ in the fibers, 
$\delta\delta^*+\delta^*\delta=(k+l)\textrm{id}$ (one can compute it in coordinates or use the Cartan formula for the Lie 
derivative by $\Theta_E$). So we define the operator $\kappa$ to be $\frac{1}{k+l}\delta^*$ on $E$-$k$-differential forms 
with value in $\mathcal B$ and $l$-polynomial in the fibers for $k+l>0$, and $0$ on sections of $\mathcal B$ constant in 
the fibers. Then one has 
\begin{equation}\label{eq:homotopy}
u=\delta\kappa u+\kappa\delta u+\mathcal Hu\qquad u\in~^E\Omega(X,\mathcal B)
\end{equation}
where $\mathcal Hu=u_{\vert y^i=\xi^i=0}\in F^0\mathcal B$ is the \emph{harmonic} part of $u$. 
\end{proof}

\subsection{Flattening the connection}

Let $\nabla$ be a linear torsion free $E$-connection (it always exists). 

The connection defines a derivation of $^E\Omega(X,\mathcal W)$ (which we will identify by the same symbol $\nabla$). 
Denote by $\Gamma_{ij}^k(x)$ Christoffel's symbols of $\nabla$; thus one can write the induced derivation in local 
coordinates 
$$\nabla=d_E+\Gamma$$
where $d_E$ is as in section \ref{sec:pre} and $\Gamma=-\xi^i\Gamma_{ij}^k(x)y^j\frac\partial{\partial y^k}$. 

This derivation $\nabla$ obviously extends to derivations of the DGLA $^E\Omega(X,\mathcal T)$ and 
$^E\Omega(X,\mathcal D)$. Namely 
$$\nabla=d_E+[\Gamma,\cdot]_S:~^E\Omega^*(X,\mathcal T)\to~^E\Omega^{*+1}(X,\mathcal T)$$
$$\nabla=d_E+[\Gamma,\cdot]_G:~^E\Omega^*(X,\mathcal D)\to~^E\Omega^{*+1}(X,\mathcal D)$$
On one hand it is clear by definition that $\nabla$ is indeed a derivation of the Lie super-algebra structures. On the 
other hand $d_E(m)=0$ and $[\Gamma,m]_G=0$ (this is just Leibniz rule), and hence $\nabla$ super-commutes with 
$\partial$. 

Since the connection is torsion free (i.e., $\Gamma_{ij}^k-\Gamma_{ji}^k=c_{ij}^k$), $\nabla$ and $\delta$ super-commute: 
$$\nabla\delta+\delta\nabla=\xi^i\wedge\xi^j(\Gamma_{ij}^k(x)-\frac12c_{ij}^k(x))\frac\partial{\partial y^k}=0$$

The standard curvature $E$-$(1,3)$-tensor of the connection induces an operator $R$ on $^E\Omega(X,\mathcal W)$ which is 
given in local coordinates by 
$$R=-\frac12\xi^i\wedge\xi^jR_{ijk}^l(x)y^k\frac\partial{\partial y^l}:
~^E\Omega^*(X,\mathcal W)\to~^E\Omega^{*+2}(X,\mathcal W)$$
where $R_{ijk}^l$ are the coefficients of the curvature $E$-tensor (\ref{eq:conn}). 
Then one has 
$$\nabla^2=\xi^i\wedge\xi^j(\Gamma_{ik}^m\Gamma_{jm}^l+\rho(e_j)\cdot\Gamma_{ik}^l+\frac12c_{ij}^m\Gamma_{mk}^l)
y^k\frac\partial{\partial y^l}=R$$
Obviously $\nabla^2$ acts as $[R,\cdot]_S$ and $[R,\cdot]_G$ respectively on $^E\Omega(X,\mathcal T)$ and 
$^E\Omega(X,\mathcal D)$. \\

Eventhough $\nabla$ is not square zero in general, we use it to deform the differential $\delta$. Namely, using an element 
\begin{equation}\label{eq:flat}
A=\sum_{p=2}^\infty\xi^kA_{ki_1\dots i_p}^jy^{i_1}\cdots y^{i_p}\frac\partial{\partial y^j}
\in~^E\Omega^1(X,\mathcal T^0)\subset~^E\Omega^1(X,\mathcal D^0)
\end{equation}
we construct a new derivation 
\begin{eqnarray}
D=\nabla-\delta+A:~^E\Omega^*(X,\mathcal W)\to~^E\Omega^{*+1}(X,\mathcal W) && \nonumber \\
D=\nabla-\delta+[A,\cdot]_S:~^E\Omega^*(X,\mathcal T)\to~^E\Omega^{*+1}(X,\mathcal T) && \label{eq:der} \\
D=\nabla-\delta+[A,\cdot]_G:~^E\Omega^*(X,\mathcal D)\to~^E\Omega^{*+1}(X,\mathcal D) && \nonumber
\end{eqnarray}

In some sens, $\nabla$ can be viewed as a connection on the "big" bundles $\mathcal B$ which we flatten recursively 
by adding terms of higher polynomial degree in the fibers. 
\begin{prop}\label{thm:flat}
There exists an element $A$ as in (\ref{eq:flat}) such that $\kappa A=0$ and the corresponding derivation $D$ 
(\ref{eq:der}) is square zero, $D^2=0$. 
\end{prop}
In what follow, we write $[A,\cdot]$ for $A\cdot$, $[A,\cdot]_S$, $[A,\cdot]_G$ when $\mathcal B$ is respectively 
$\mathcal W$, $\mathcal T$, $\mathcal D$. 
\begin{proof}
Since $\kappa$ raises the polynomial degree in the fibers (i.e., in $y$), there is a unique solution $A$ in the form 
(\ref{eq:flat}) to equation 
\begin{equation}\label{eq:A1}
A=\kappa R+\kappa(\nabla A+\frac12[A,A])
\end{equation}
First observe that $\kappa^2=0$ implies that $\kappa A=0$. 

Now let us show that $A$ satisfies equation
\begin{equation}\label{eq:A2}
\delta A=R+\nabla A+\frac12[A,A]
\end{equation}
which obviously implies that $D^2=0$. \\
Using (\ref{eq:homotopy}) together with $\kappa A=0=\mathcal HA$ one finds that 
\begin{equation}\label{eq:A3}
\kappa\delta A=\kappa R+\kappa(\nabla A+\frac12[A,A])
\end{equation}
Define $C=-\delta A+R+\nabla A+\frac12[A,A]$. Due to (\ref{eq:A3}) $\kappa C=0$, and reformulating Bianchi's identities 
for $R$ one can show that $\delta R=0=\nabla R$. These equalities, together with (\ref{eq:homotopy}), imply that 
$C=\kappa(\nabla C+[A,C])$. Since the operator $\kappa$ raises the polynomial degree in the fiber, this latter 
equation has a unique zero solution. \\
Thus $A$ satisfies (\ref{eq:A2}) and the proposition is proved. 
\end{proof}

\subsection{Acyclicity of the complexes. Resolutions}\label{sec:reso}

\begin{thm}\label{thm:acyclic}
$H^*(^E\Omega(X,\mathcal B),D)=H^0(^E\Omega(X,\mathcal B),D)\cong F^0\mathcal B$. 
\end{thm}
\begin{proof}
Using arguments similar to those of the proof of proposition \ref{thm:flat}, one can show that 
$H^*(^E\Omega(X,\mathcal B),D)=H^0(^E\Omega(X,\mathcal B),D)$ (see \cite[theorem 3]{D} for details). 

Let us now prove that $H^0(^E\Omega(X,\mathcal B),D)\cong F^0\mathcal B$. For any $u_0\in F^0\mathcal B$, there is a 
unique solution $u\in~^E\Omega^0(X,\mathcal B)=\Gamma(X,\mathcal B)$ of the equation
\begin{equation}\label{eq:u1}
u=u_0+\kappa(\nabla u+[A,u])
\end{equation}
(still because $\kappa$ raises the polynomial degree in the fibers). It is obvious that $\mathcal Hu=u_0$; let us 
prove that $Du=0$. 
Let $v=Du$, then $Dv=0=\mathcal Hv$, and $\kappa v=0$ after (\ref{eq:u1}). Then use (\ref{eq:homotopy}) to find 
$$v=\kappa(\nabla v+[A,v])$$
Again, this equation has a unique zero solution and consequently $v=Du=0$. 

We have defined a linear map $\vartheta:F^0\mathcal B\to Z^0(^E\Omega(X,\mathcal B),D)=H^0(^E\Omega(X,\mathcal B),D)$ that 
sends $u_0$ to the solution $u=\vartheta(u_0)$ of (\ref{eq:u1}) and such that $\mathcal H(\vartheta(u_0))=u_0$. 
This map is obviously injective: the solution of (\ref{eq:u1}) is zero if and only if $u_0=0$. It is also surjective: 
if $v\in Z^0(^E\Omega(X,\mathcal B),D)$ is such that $\mathcal Hv=0$ then by (\ref{eq:homotopy}) 
$v=\kappa(\nabla v+[A,v])$ and so $v=0$. 
\end{proof}

In the case of the Weyl bundle $\mathcal B=\mathcal W$, one can easily show that 
$$\mathcal H:H^*(^E\Omega(X,\mathcal W),D)=Z^0(^E\Omega(X,\mathcal W),D)\to F^0\mathcal W=C^\infty(X)$$
is a morphism of commutative algebras. 

In the same spirit we have $F^0\mathcal T=\Gamma(X,\wedge E)=T_{poly}E$. On the other hand the differential $D$ respects 
the DGLA structure on $^E\Omega(X,\mathcal T)$ and thus its homology acquires a DGLA structure. In the following 
proposition we show that $\mathcal H$ respects the DGLA structures. 
\begin{prop}\label{thm:reso}
$H^*(^E\Omega(X,\mathcal T),D))\cong_{DGLA}T_{poly}E$
\end{prop}
\begin{proof}
Let $u,v\in Z^0(^E\Omega(X,\mathcal T),D)$. We are going to show that 
\begin{equation}\label{eq:Lie}
\mathcal H([u,v]_S)=[\mathcal H(u),\mathcal H(v)]_E
\end{equation}
Since $\mathcal H$ preserves the exterior product of polyvector fields it is sufficient to prove it in the following two 
cases: first when $u_0=\mathcal H(u)$ and $v_0=\mathcal H(v)$ are vector fields, next when $u_0=\mathcal H(u)$ a vector 
field and $f=\mathcal H(v)$ is a function. \\
\emph{First case. }Let $u_0=u^i(x)\frac\partial{\partial y^i}$ and $v_0=v^i(x)\frac\partial{\partial y^i}$ be vector 
fields. Note that 
$$[u_0,v_0]_E=(u^i\rho(e_i)v^k+u^iv^jc_{ij}^k-v^i\rho(e_i)u^k)\frac\partial{\partial y^k}$$
Then, by an easy calculation one obtains 
$$\begin{array}{rccccl}
u&=&\vartheta(u_0)&=&u+y^i(\rho(e_i)u^k+\Gamma_{ij}^ku^j)\frac\partial{\partial y^k} & \textrm{mod}~\vert y\vert^2 \\
v&=&\vartheta(v_0)&=&v+y^i(\rho(e_i)v^k+\Gamma_{ij}^kv^j)\frac\partial{\partial y^k} & \textrm{mod}~\vert y\vert^2 \\
\end{array}$$
And thus 
$$\begin{array}{rcll}
[u,v]_S & = & 
u^i(\rho(e_i)v^k+\Gamma_{ij}^kv^j)\frac\partial{\partial y^k}-v^i(\rho(e_i)u^k+\Gamma_{ij}^ku^j)\frac\partial{\partial y^k}
& \textrm{mod}~\vert y\vert \\
& = & (u^i\rho(e_i)v^k+c_{ij}^ku^iv^j-v^i\rho(e_i)u^k)\frac\partial{\partial y^k} & \textrm{mod}~\vert y\vert \\
& = & \vartheta([u_0,v_0]_E) & \textrm{mod}~\vert y\vert
\end{array}$$
\emph{Second case. }Let $u_0=u^i(x)\frac\partial{\partial y^i}$ be a vector field and $f$ be a function. One has 
$$[u_0,f]_E=\rho(u_0)f=u^i\rho(e_i)f$$
Since $v=\vartheta(f)=f+y^i\rho(e_i)f~\textrm{mod}~\vert y\vert^2$ we obtain 
$$[u,v]_S=u^i\rho(e_i)f~\textrm{mod}~\vert y\vert$$
Consequently (\ref{eq:Lie}) is satisfied and $\mathcal H$ is an isomorphism of graded Lie algebras. Since the 
differentials are both zero it is a DGLA-isomorphism. 
\end{proof}
As above $D$ preserves the DGLA structure on $^E\Omega(X,\mathcal D)$ and thus its homology is also a DGLA. Using the 
PBW theorem for Lie algebroids (see \cite{R,NWX}) one finds that $F^0\mathcal D=\Gamma(X,\otimes^*S(E))$ and 
$D_{poly}E=\widetilde\otimes^*U(E)$ are isomorphic as (filtered) vector spaces. Again we have: 
\begin{prop}\label{thm:reso2}
$H^*(^E\Omega(X,\mathcal D),D)\cong_{DGLA}D_{poly}E$
\end{prop}
\begin{proof}
Let us first set $\tau_0=1$ and $\tau_{k+1}=y^ie_i\tau_k-y^iy^j\Gamma_{ij}^l(x)\frac{\partial\tau_k}{\partial y^l}$. For any 
$k\in\N$ $\tau_k$ is a well-defined element of $\Gamma(X,\mathcal W)\otimes_{C^\infty(X)}UE$. Then for any fiberwise 
differential operator constant in the fibers 
$u=u^{i_1\dots i_k}(x)\frac\partial{\partial y^{i_1}}\cdots\frac\partial{\partial y^{i_k}}$ we define 
$\mu(u):=\frac1{k!}(\mathcal H\otimes\textrm{id})(u\cdot\tau_k)$ and compute 
$$\mu(u)=\frac1{k!}u^{i_1\dots i_k}(x)\sum_{\sigma\in S_k}e_{i_{\sigma_1}}\cdots e_{i_{\sigma_k}}\quad\textrm{mod}~U_{k-1}(E)$$
Thus $\mu$ is filtered and its associated graded map coincides with the usual isomorphism 
$\Gamma(X,S^k(E))\tilde\to U_k(E)/U_{k-1}(E)$ (see \cite[proof of theorem 3]{NWX}). Consequently $\mu$ is an isomorphism 
which naturaly extends to an isomorphism from $F^0\mathcal D$ to $D_{poly}E$. Moreover, $\mu\circ\Delta=\Delta\circ\mu$ on 
$\Gamma(X,S(E))$. 

Next we also have $\mathcal H\circ\Delta=\Delta\circ\mathcal H$ on $Z^0(^E\Omega(X,\mathcal D^0),D)$ and thus the 
composition $\mu\circ\mathcal H:Z^0(^E\Omega(X,\mathcal D),D)\to D_{poly}E$ is compatible with coproducts. In particular 
it commutes with differentials: $(\mu\circ\mathcal H)\circ\partial=\partial\circ(\mu\circ\mathcal H)$. 

Finally, due to the compatibility with coproducts and to the special form of the brackets (see \ref{sec:pre} and 
\ref{sec:A2}) it is now sufficient to show that $\mu\circ\mathcal H$ restricts to a morphism of associative algebras 
between $Z^0(^E\Omega(X,\mathcal D^0),D)$ and $U(E)$; we have to prove that 
$P_1P_2=\mu\circ\mathcal H(\widetilde P_1\widetilde P_2)$ where $P_i$ are generators of $U(E)$ and 
$\widetilde P_i=\vartheta\circ\mu^{-1}(P_i)$. 
There are four distinct cases: \\
\emph{First case. }$P_1=u^i(x)e_i$ and $P_2=v^j(x)e_j$ are vector fields. Then 
$\mu^{-1}(P_1)=\nolinebreak u^i\frac\partial{\partial y^i}$, $\mu^{-1}(P_2)=v^j\frac\partial{\partial y^j}$ and thus 
$\mathcal H(\widetilde P_1\widetilde P_2)
=u^iv^j\frac{\partial^2}{\partial y^i\partial y^j}+u^i(\rho(e_i)v^k+\Gamma_{ij}^kv^j)\frac\partial{\partial y^k}$. 
And since $\mu(u^iv^j\frac{\partial^2}{\partial y^i\partial y^j})
=\frac12u^iv^j(e_ie_j+e_je_i-(\Gamma_{ij}^k+\Gamma_{ji}^k)e_k)$ one computes 
\begin{eqnarray*}
P_1P_2 & = & u^ie_iv^je_j=\frac12(u^iv^j(e_ie_j+e_je_i+c_{ij}^ke_k))+u^i\rho(e_i)v^je_j \\
& = & \mu(u^iv^j\frac{\partial^2}{\partial y^i\partial y^j}+u^i(v^j\Gamma_{ij}^k+\rho(e_i)v^k)\frac\partial{\partial y^k})
=\mu\circ\mathcal H(\widetilde P_1\widetilde P_2)
\end{eqnarray*}
\emph{Second and third cases. }$P_1=u^i(x)e_i$ is a vector field and $P_2=f$ is a function. Since $\mu^{-1}(P_1)=f$ we have 
$P_1P_2=u^i(fe_i+\rho(e_i)f)=\mu(u^i(f\frac\partial{\partial y^i}+\rho(e_i)f))=
\mu\circ\mathcal H(\widetilde P_1\widetilde P_2)$ and 
$P_2P_1=fu^ie_i=\mu(fu^i\frac\partial{\partial y^i})=\mu\circ\mathcal H(\widetilde P_2\widetilde P_1)$. \\
\emph{Fourth case. }$P_1=f$ and $P_2=g$ are functions. $P_1P_2=fg=\mu\circ\mathcal H(\widetilde P_1\widetilde P_2)$. \\
Consequently $\mu\circ\mathcal H$ is a DGLA-isomorphism. 
\end{proof}

\section{Proof of theorem \ref{thm:formal}}

\subsection{Twisting a fiberwise quasi-isomorphism}\label{sec:twist}

In virtue of properties 1 and 2 in theorem \ref{thm:K} we have a fiberwise quasi-isomorphism $U_K$ from 
$(^E\Omega(X,\mathcal T),0,[,]_S)$ to $(^E\Omega(X,\mathcal D),\partial,[,]_G)$. Our purpose is to twist $U_K$ 
in order to get a quasi-isomorphism $\widetilde U$ from $(^E\Omega(X,\mathcal T),D,[,]_S)$ to 
$(^E\Omega(X,\mathcal T),\partial+D,[,]_G)$. 

Let us recall that the differential $D$ can be written locally in the form 
\begin{eqnarray*}
& D=d_E+[B,\cdot]_S:~^E\Omega^*(X,\mathcal T)\to~^E\Omega^{*+1}(X,\mathcal T) & \\
& D=d_E+[B,\cdot]_G:~^E\Omega^*(X,\mathcal D)\to~^E\Omega^{*+1}(X,\mathcal D) & \\
& \textrm{where }B=-\xi^i\frac\partial{\partial y^i}-\xi^i\Gamma_{ij}^k(x)y^j\frac\partial{\partial y^k}+
\sum_{p\geq2}\xi^iA_{ij_1\dots j_p}^k(x)y^{j_1}\cdots y^{j_p}\frac\partial{\partial y^k} & 
\end{eqnarray*}
Let $V$ be a $E$-coordinates disk, then we prove 
\begin{prop}
$U_K$ defines a quasi-isomorphism of DGLA from $(^E\Omega(V,\mathcal T),d_E,[,]_S)$ to 
$(^E\Omega(V,\mathcal D),\partial+d_E,[,]_G)$. 
\end{prop}
\begin{proof}
Let us note respectively $\mathbb T$ and $\mathbb D$ for $^E\Omega(V,\mathcal T)$ and $^E\Omega(V,\mathcal D)$. 
Since $d_E$ commutes with the fiberwise DGLA structures of $\mathbb T$ and $\mathbb D$, and also with the fiberwise 
$L_\infty$-morphism $U_K$, then $U_K$ defines a $L_\infty$-morphism from $(\mathbb T,d_E,[,]_S)$ to 
$(\mathbb D,\partial+d_E,[,]_G)$. 

Now observe that $(\mathbb T,0,d_E)$ and $(\mathbb D,\partial,d_E)$ are double complexes; $U_K^{[1]}=A_{hkr}$ is an 
inclusion of double complexes. Thus we have a long exact sequence in cohomology 
$$\cdots\to H^k(\mathbb T,d_E)\to H^k(\mathbb D,\partial+d_E)\to H^k(\mathbb D/\mathbb T,\partial+d_E)\to\cdots$$
Since the inclusion $A_{hkr}:(\mathbb T,0)\to(\mathbb D,\partial)$ is a quasi-isomorphism of complexes one has 
$H^*(\mathbb D/\mathbb T,\partial)=0$. Then the (second) spectral sequence of the double complex 
$(\mathbb D/\mathbb T,\partial,d_E)$ goes to zero and thus $H^*(\mathbb D/\mathbb T,\partial+d_E)=0$. 

Consequently $A_{hkr}$ induces an isomorphism $H^*(\mathbb T,d_E)\tilde\rightarrow H^*(\mathbb D,\partial+d_E)$. It means 
that $U_K$ is a quasi-isomorphism of DGLA from $(\mathbb T,d_E,[,]_S)$ to $(\mathbb D,\partial+d_E,[,]_G)$. 
\end{proof}
On $V$ the element $B\in~\!\!\!^E\Omega(V,\mathcal T^0)\subset~\!\!\!^E\Omega(V,\mathcal D^0)$ is 
well-defined, and since $D$ is square zero it is a Maurer-Cartan element. It means that $(^E\Omega(V,\mathcal T),D,[,]_S)$ 
(resp.~$(^E\Omega(V,\mathcal D),\partial+D,[,]_G)$) can be obtained using a twisting of the DGLA 
$(^E\Omega(V,\mathcal T),d_E,[,]_S)$ (resp.~$(^E\Omega(V,\mathcal D),\partial+d_E,[,]_G)$) by $B$ (a general description 
of twisting procedures for DGLA and their $L_\infty$-morphisms is presented in \cite[section 2.3]{D2}). \\
Due to properties 3 and 5 in theorem \ref{thm:K}, $U_K$ maps $B\in~\!\!\!^E\Omega(V,\mathcal T^0)$ to 
$B\in~\!\!\!^E\Omega(V,\mathcal D^0)$. Then we can define a quasi-isomorphism $\widetilde U$ of DGLA from 
$(^E\Omega(V,\mathcal T),D,[,]_S)$ to $(^E\Omega(V,\mathcal T),\partial+D,[,]_G)$ using a twisting of $U_K$ by $B$. 
Namely, 
\begin{equation}\label{eq:twisting}
\widetilde U(Y)=\textrm{exp}((-B)\wedge)U_K(\textrm{exp}(B\wedge)Y)
\end{equation}
Next proposition tells us that $\widetilde U$ does not depend on the choice of local coordinates and thus is defined 
globaly. 
\begin{prop}
$\widetilde U$ extends to a quasi-isomorphism of DGLA from $(^E\Omega(X,\mathcal T),D,[,]_S)$ to 
$(^E\Omega(X,\mathcal D),\partial+D,[,]_G)$. 
\end{prop}
\begin{proof}
See \cite[proposition 3]{D} (it makes use of property 4 in theorem \ref{thm:K}). 
\end{proof}

\subsection{End of the proof: contraction of $\widetilde U$}\label{sec:contract}

On one hand we know from proposition \ref{thm:reso} that there exists a quasi-isomorphism of DGLA $U_T$ from 
$(T_{poly}E,0,[,]_E)$ to $(^E\Omega(X,\mathcal T),D,[,]_S)$. On the other hand we have also a quasi-isomorphism of DGLA 
$\widetilde U$ from $(^E\Omega(X,\mathcal T),D,[,]_S)$ to $(^E\Omega(X,\mathcal D),\partial+D,[,]_G)$ (section 
\ref{sec:twist}). Let us define $\overline U=\widetilde U\circ U_T$ and claim
\begin{prop}
One can modify $\overline U$ to construct a quasi-isomorphism of DGLA $\underline U$ from $T_{poly}E$ to 
$^E\Omega(X,\mathcal D)$ whose structure maps take values in $Z^0(^E\Omega(X,\mathcal D),D)$. 
\end{prop}
\begin{proof}
See \cite[proposition 5]{D}. 
\end{proof}
Consequently, composing $\underline U$ with the DGLA-isomorphism of proposition \ref{thm:reso2} we obtain a 
quasi-isomorphism of DGLA $U_E$ from $(T_{poly}E,0,[,]_E)$ to $(D_{poly}E,\partial,[,])$. Thus we have proved theorem 
\ref{thm:formal}. Q.E.D.

\subsection{Equivariant formality theorem}

By a \emph{good} action of a group $G$ on a Lie algebroid $E\to X$ we mean a smooth action, linear in the fibers and 
compatible with anchor map and bracket: for all $g\in G$, $u,v\in E$, $g\cdot\rho(u)=\rho(g\cdot u)$ and 
$g\cdot[u,v]_E=[g\cdot u,g\cdot v]_E$. Such an action extends naturally to $T_{poly}E$ and $D_{poly}E$ with the property 
that all structures are $G$-invariant. 
In this context the quasi-isomorphism of complexes defined in theorem \ref{thm:hkr} is $G$-equivariant. In particular 
it restricts to a quasi-isomorphism of complexes $U_{hkr}:(T_{poly}E)^G\to(D_{poly}E)^G$. 
The following theorem is a $G$-equivariant version of theorem \ref{thm:formal}. 
\begin{thm}\label{thm:equiv}
Consider a Lie algebroid $E\to X$ equipped with a good action of a group $G$. If there exists a $G$-invariant torsion free 
$E$-connexion $\nabla$, then one can construct a $G$-equivariant quasi-isomorphism of DGLA from $T_{poly}E$ to $D_{poly}E$. 
\end{thm}
\begin{proof}
First one can canonically extend the action of $G$ to the spaces $^E\Omega(M,\mathcal W)$, $^E\Omega(M,\mathcal T)$ and 
$^E\Omega(M,\mathcal D)$ in such a way that all algebraic structures we have defined are $G$-invariant. 

First we are going to prove that the resolutions constructed in section \ref{sec:reso} are $G$-equivariant. 
The differential $\delta$, the homotopy operator $\kappa$ and the projection $\mathcal H$ are obviously $G$-invariant. 
The $G$-invariance of the connection $\nabla$ implies the $G$-invariance of the induced derivation (also called $\nabla$) 
and of the curvature tensor $R$. Thus equation (\ref{eq:A1}) has a $G$-invariant solution $A$ (\ref{eq:flat}) and then 
the differential $D$ (\ref{eq:der}) of proposition \ref{thm:flat} is $G$-invariant. In the same way, $\vartheta$ is 
$G$-invariant since it is defined by $G$-equivariant equation (\ref{eq:u1}). Thus the DGLA-isomorphisms of propositions 
\ref{thm:reso} and \ref{thm:reso2} are $G$-equivariant. 

Second, since $G$ acts on the fibers by linear transformations and due to property 2 in theorem \ref{thm:K} the fiberwise 
quasi-isomorphism $U_K$ is $G$-equivariant. 

Third we have to prove that the quasi-isomorphim $\widetilde U$ constructed with the help of the twisting procedure 
(\ref{eq:twisting}) is $G$-equivariant. Let $V$ be a coordinate disk and $B$ be the twisting element of section 
\ref{sec:twist}. Since $U_K$ is $G$-equivariant one has 
$$gU_K^{[n+m]}(B,\dots,B,v_1,\dots,v_n)=U_K^{[n+m]}(gB,\dots,gB,gv_1,\dots,gv_n)$$
where $g\in G$ acts on $B$ as on a tensor element. Now a sufficient condition for $\widetilde U$ to be $G$-equivariant is 
$$gU_K^{[n+m]}(B,\dots,B,v_1,\dots,v_n)=U_K^{[n+m]}(g\cdot B,\dots,g\cdot B,gv_1,\dots,gv_n)$$
where $g\cdot$ acts by ususal transformations of Christoffel's symbols in $B$. 
Then remark that $g\cdot B-gB$ is a fiberwise polyvector field linear in the fibers on $V$; thus using property 4 of 
theorem \ref{thm:K} we obtain the desired result. 

Finally, it is not difficult to see that the contraction procedure of section \ref{sec:contract} involves only 
$G$-equivariant cohomological equations (see \cite{D} for details). 
\end{proof}
\begin{exs}\emph{
(i) Consider the case of a Lie algebra $\al$ (i.e., a Lie algebroid over a point) with the adjoint action of its Lie group 
$G$ (which is a good action). Then the Lie algebroid connection given by half the Lie bracket on $\al$ is a 
torsion free $G$-invariant connection and we obtain a $G$-equivariant quasi-isomorphism of DGLA from $\wedge^*\al$ to 
$\otimes^*U\al$. In particular for any subgroup $H\subset G$ one obtains a quasi-isomorphism of DGLA from $(\wedge^*\al)^H$ 
to $(\otimes^*U\al)^H$. }

\emph{(ii) If a group $G$ acts smoothly on a manifold $X$, then it induces a good action on $TX$. In this particular case 
our theorem is equivalent to theorem 5 of \cite{D}. }

\emph{(iii) Now if $E\to X$ is a Lie algebroid with injective anchor (i.e., $E$ is the Lie algebroid of a foliation), then 
any smooth action of a group $G$ on $X$ that respects the foliation (i.e., that sends a leaf to a leaf) gives rise to a 
good action on $E$. In this context we obtain a leafwise version of the previous example. }
\end{exs}

\appendix

\section{Formality, $L_\infty$ and all that}

\subsection{Quasi-isomorphisms of differential graded Lie algebras}

Let $(\al,d,[,])$ be a differential graded Lie algebra (DGLA). We assume that the differential is of 
degree one and the Lie super-bracket is of degree zero. One can associate to $\al$ a cocommutative coalgebra 
$C_*(\al[1])$ cofreely generated by the vector space $\al$ with a shifted parity, equipped with a coderivation $Q$ 
having two non-vanishing structure maps $Q^{[1]}=d:\al\to\al[1]$ and $Q^{[2]}=[,]:\wedge^2\al\to\al$. The fact that 
$(\al,d,[,])$ is a DGLA is equivalent to the nilpotency of $Q$ (i.e., $Q^2=0$). 
\begin{dfn}\emph{
A \emph{$L_\infty$-morphism} between two DGLA $(\gun,d_1,[,]_1)$ and $(\gde,d_2,[,]_2)$ is a morphism of 
cocommutative coalgebras $L:C_*(\gun)\to C_*(\gde)$ compatible with the DGLA structures in the following sens: 
$Q_2\circ L=L\circ Q_1$, where $Q_i$ is the square zero coderivation corresponding to $(d_i,[,]_i)$. 
}\end{dfn}
\begin{dfn}\emph{
A \emph{quasi-isomorphism of DGLA} from $(\gun,d_1,[,]_1)$ to $(\gde,d_2,[,]_2)$ is a $L_\infty$-morphism $U$ 
from $\gun$ to $\gde$ whose first structure map $U^{[1]}:\gun\to\gde$ induces an isomorphism in cohomology 
$H^*(\gun,d_1)\cong H^*(\gde,d_2)$. 
}\end{dfn}
A DGLA is \emph{formal} if it is quasi-isomorphic to the graded Lie algebra (with zero differential) of its cohomology. 

\subsection{Kontsevich formality theorem}\label{sec:A2}

Let $D_{poly}(X)$ be the vector space of polydifferential operators on a smooth manifold $X$. It is a graded vector space 
$$D_{poly}(X)=\bigoplus_{k\geq-1}D^k_{poly}(X)$$
where $D^k_{poly}(X)$ denotes the subspace of operators of rank $k+1$. 
We define on $D_{poly}(X)$ a Lie super-bracket (the Gerstenhaber bracket) given on homogeneous 
elements by $P_i\in D^{k_i}_{poly}(X)~(i=1,2)$ by $[P_1,P_2]_G=P_1\bullet P_2-(-1)^{k_1k_2}P_2\bullet P_1$, where 
$$P_1\bullet P_2(f_0,\dots,f_{k_1+k_2})=\sum_{i=0}^{k_1}(-1)^{ik_2}P_1(f_0,\dots,f_{i-1},
P_2(f_i,\dots,f_{i+k_2}),\dots,f_{k_1+k_2})$$
Associativity condition for the multiplication operator $m_0\in D^1_{poly}(X)$ can be written in terms of the Gerstenhaber 
bracket as $[m_0,m_0]_G=0$. Thus $(\partial=[m_0,\cdot]_G,[,]_G)$ defines a DGLA structure on $D_{poly}(X)$. 

Let now $T_{poly}(X)$ be the DGLA of polyvector fields on $X$: 
$$T_{poly}(X)=\bigoplus_{k\geq-1}T^k_{poly}(X)=\bigoplus_{k\geq-1}\Gamma(X,\wedge^{k+1}TX)$$
endowed with the standard Schouten-Nijenhuis bracket and zero differential. 

Hochschild-Kostant-Rosenberg theorem says that the antisymmetrisation map $A_{hkr}:T_{poly}(X)\to D_{poly}(X)$ induces an 
isomorphism $H^*(D_{poly}(X),\partial)\cong T_{poly}^*(X)$, and Kontsevich has proved in 
\cite{K} that $D_{poly}(X)$ is formal. We will use a version of this result when $X=\R^d_{formal}$: 
\begin{thm}[Kontsevich,\cite{K}]\label{thm:K}
There exists a quasi-isomorphism of DGLA $U_K$ from $T_{poly}(\R^d)$ to $D_{poly}(\R^d)$ which has the 
following properties: 
\begin{enumerate}
\item $U_K$ can be defined for $\R^d_{formal}$ (the formal completion of $\R^d$ at the origin) as well. 
\item $U_K$ is $GL_d(\R)$-equivariant. 
\item For any $n\geq2$, $v_1,\dots,v_n\in T^0_{poly}(\R^d_{formal})$, $U_K^{[n]}(v_1,\dots,v_n)=0$. 
\item For any $n\geq2$, $v\in\gl_d(\R)\subset T^0_{poly}(\R^d_{formal})$, $\chi_2,\dots,\chi_n\in T_{poly}(\R^d_{formal})$, 
$U_K^{[n]}(v,\chi_2,\dots,\chi_n)=0$. 
\item $U_K^{[1]}=A_{hkr}$. 
\end{enumerate}
\end{thm}

\section{Hopf algebroids}\label{sec:B}

\begin{dfn}[\cite{X2}, see also \cite{L}]\emph{
A \emph{Hopf algebroid} is an associative algebra with unit $H$ together with a base algebra $R$, an algebra homomorphism 
$s:R\to H$ and an algebra antihomomorphism $t:R\to H$ whose respective images commute together (the source and target maps, 
which give $H$ an $R$-bimodule structure), and $R$-bimodule maps $\Delta:H\to H\otimes_RH$ (the coproduct) and 
$\varepsilon:H\to R$ (the counit) such that 
\begin{enumerate}
\item $\Delta(1)=1\otimes_R1$ and $(\Delta\otimes_R\textrm{id})\circ\Delta=(\textrm{id}\otimes_R\Delta)\circ\Delta$
\item $\forall a\in R, \forall h\in H, \Delta(h)(t(a)\otimes_R1-1\otimes_Rs(a))=0$
\item $\forall h_1,h_2\in H, \Delta(h_1h_2)=\Delta(h_1)\Delta(h_2)$
\item $\varepsilon(1_H)=1_R$ and 
$(\varepsilon\otimes_R\textrm{id}_H)\circ\Delta=(\textrm{id}_H\otimes_R\varepsilon)\circ\Delta=\textrm{id}_H$
\end{enumerate}
}\end{dfn}
Given a Hopf algebroid $H$ over a base $R$, an \emph{anchor} is a representation $\rho:H\to\End R$ which is also a 
$R$-bimodule map and satisfies 
$$\begin{array}{lc}
s(\rho(x_1)\cdot a)x_2=xs(a) & x\in H,a\in R \\
x_1t(\rho(x_2)\cdot a)=xt(a) & x\in H,a\in R \\
\rho(x)\cdot1_R=\varepsilon(x) & x\in H
\end{array}$$

A \emph{twistor} (\cite{X2}) in a Hopf algebroid $H$ over a base $R$ is an invertible element $J\in H\otimes_RH$ that 
satisfies
\begin{equation}\label{eq:twist}
\begin{array}{l}
J^{12,3}J^{1,2}=J^{1,23}J^{2,3} \\
(\varepsilon\otimes_R\textrm{id})(J)=(\textrm{id}\otimes_R\varepsilon)(J)=1_H
\end{array}
\end{equation}
Let $H$ be a Hopf algebroid over a base $R$ (resp. with anchor $\rho$), and let $J=\sum_ix_i\otimes_Ry_i$ be a twistor. 
Then one can define a new product on $R$ given by $a*_Jb=\sum_i(\rho(x_i)a)(\rho(y_i)b)$, a new coproduct 
$\Delta_J=J^{-1}\Delta J$, and new source and target maps given by $s_J(a)=\sum_is(\rho(x_i)a)y_i$ and 
$t_J(a)=t(\rho(y_i)a)x_i$. Denote $R_J=(R,*_J)$. 
\begin{thm}[\cite{X2}, theorem 4.14]
Let $(H,R,\Delta,s,t,\varepsilon)$ be a Hopf algebroid (resp. with anchor $\rho$). If $J$ is a twistor, then 
$(H,R_J,\Delta_J,s_J,t_J,\varepsilon)$ is again a Hopf algebroid (resp. with the same anchor $\rho$). 
\end{thm}

\section{Lie algebroid connections}

Let $(E,[,]_E,\rho)$ be a Lie algebroid over a smooth manifold $X$. 
\begin{dfn}\emph{
A \emph{linear $E$-connection} is a map $\nabla:\Gamma(X,E)\times\Gamma(X,E)\to\Gamma(X,E)$ such that 
\begin{enumerate}
\item $\nabla$ is $C^\infty(X)$-linear with respect to the first argument. 
\item $\nabla$ is $\R$-linear with respect to the second argument. 
\item for all $f\in C^\infty(X)$ and $u,v\in\Gamma(X,E)$, $\nabla_ufv=f\nabla_uv+(\rho(u)\cdot f)v$. 
\end{enumerate}
}\end{dfn}
In a local base $(e_1,\dots,e_r)$ of $E$, $\nabla$ is completely determined by its \emph{Christoffel's symbols} 
$\Gamma_{ij}^k$ which are given by: $\nabla_{e_i}e_j=\Gamma_{ij}^ke_k$. \\
\emph{Remark. }As with usual connections, one can define the covariant derivative on $E$-tensor in a unique way such that 
$\nabla_u$ is a derivation with respect to the tensor product of $E$-tensors, commutes with the contraction of $E$-tensors, 
acts as $\rho(u)$ on functions, and is $\R$-linear. 
\begin{dfn}\emph{
(i) The \emph{torsion} $T$ of $\nabla$ is the $E$-$(1,2)$-tensor defined by 
$$T(u,v)=\nabla_uv-\nabla_vu-[u,v]_E$$
\indent(ii) The \emph{curvature} $R$ of $\nabla$ is the $E$-$(1,3)$-tensor defined by 
$$R(u,v)w=([\nabla_u,\nabla_v]-\nabla_{[u,v]_E})w$$
}\end{dfn}
\noindent Coefficients of these tensors can be expressed in a local base $(e_1,\dots,e_n)$: 
\begin{equation}\label{eq:conn}
\begin{array}{l}
T_{ij}^k=\Gamma_{ij}^k-\Gamma_{ji}^k-c_{ij}^k \\
R_{ijk}^l=\Gamma_{im}^l\Gamma_{jk}^m-\Gamma_{ik}^m\Gamma_{jm}^l+\rho(e_i)\cdot\Gamma_{jk}^l
-\rho(e_j)\cdot\Gamma_{ik}^l-c_{ij}^m\Gamma_{mk}^l
\end{array}
\end{equation}
\begin{prop}
There exists a torsion free linear $E$-connection. 
\end{prop}
\begin{proof}
Let $(U_\alpha)_\alpha$ be a covering of $X$ by trivializing opens for $E$. On each $U_\alpha$ one has a basis $(e_i)_i$ of sections 
and then can define $\nabla^{(\alpha)}_{e_i}e_j=\frac12[e_i,e_j]$. Let $f_\alpha$ be such that $\sum_\alpha f_\alpha=1$ and 
define $\nabla=f_\alpha\nabla^{(\alpha)}$. $\nabla$ is a torsion free linear $E$-connection. 
\end{proof}
\begin{prop}[Bianchi's identities]
For all $u,v,w\in\Gamma(X,E)$ 
$$\nabla_uR(v,w)+R(T(u,v),w)+c.p.(u,v,w)=0$$
$$R(u,v)w-T(T(u,v),w)-\nabla_uT(v,w)+c.p.(u,v,w)=0$$
\end{prop}
\begin{proof}
See for example \cite{Fer}. 
\end{proof}

\end{document}